# Heteroclinic orbits, and transport in a perturbed Integrable Standard map

Héctor E. Lomelí and James D. Meiss *
Program in Applied Mathematics
University of Colorado
Boulder, CO 80309-0526
lomeli@boulder.colorado.edu and jdm@boulder.colorado.edu

Feb 27, 1996

*JDM acknowledges support from the NSF under grant DMS-9305847

### Abstract

Explicit formulae are given for the saddle connection for an integrable family of standard maps studied by Suris. A generalization of Melnikov's method shows that, upon perturbation, this connection is destroyed. We give explicit formula for the first order approximation of the area of the lobes of the resultant turnstile. It is shown that the lobe area is exponentially small in the limit when the Suris map approaches the trivial twist map.

## List of Figures

# 1 Introduction

Standard maps are area-preserving diffeomorphisms of $\mathbb{T} \times \mathbb{R}$ given by

$$f(\theta, r) = (\theta + r + V'(\theta), r + V'(\theta)) \ , \tag{1}$$

where the potential, $V$, is periodic, $V(\theta+1) = V(\theta)$. The case where

$$V(\theta) = \frac{k}{4\pi^2}\cos(2\pi\theta) \tag{2}$$

is known as *the* standard or Taylor-Chirikov map. This model is important because it gives a local description of nonintegrable two degree of freedom Hamiltonian dynamics.

Twist maps, of which the standard map is an example, will be our major concern in this paper (for review, see [17]). Such maps have Lagrangian generating functions, $S(\theta, \theta')$, which generate the map implicitly through the equations

$$\begin{aligned} r = & \ -\partial_1 S(\theta, \theta') \ , \\ r' = & \ \ \partial_2 S(\theta, \theta') \ . \end{aligned} \tag{3}$$

To generate a map, $S$ must satisfy the *twist condition* that the second equation above can be inverted to obtain $\theta'(r, \theta)$; this occurs, e.g. if

$$\partial_1 \partial_2 S < 0 \ , \tag{4}$$

and implies the geometric condition that vertical lines tilt to the right upon iteration, $\frac{\partial \theta'}{\partial r} > 0$. Furthermore, we assume that our twist map has *zero net flux*, which is equivalent to

$$S(\theta+1, \theta'+1) = S(\theta, \theta') \ . \tag{5}$$

For maps of the standard form, the generating function is

$$S(\theta, \theta') = \frac{1}{2}(\theta' - \theta)^2 + V(\theta) \ . \tag{6}$$

In this paper we discuss a standard map introduced by Suris [21, 17] with the potential given by

$$\begin{aligned} V(\theta) = & \ -\tfrac{2}{\pi}\int_0^\theta dt \, \tan^{-1}\left(\tfrac{\delta \sin(2\pi t)}{1+\delta\cos(2\pi t)}\right) \\ = & \ \tfrac{1}{\pi^2}\Re\left[\mathrm{dilog}(1+\delta) - \mathrm{dilog}(1+\delta e^{2\pi i\theta})\right] \ , \end{aligned}$$

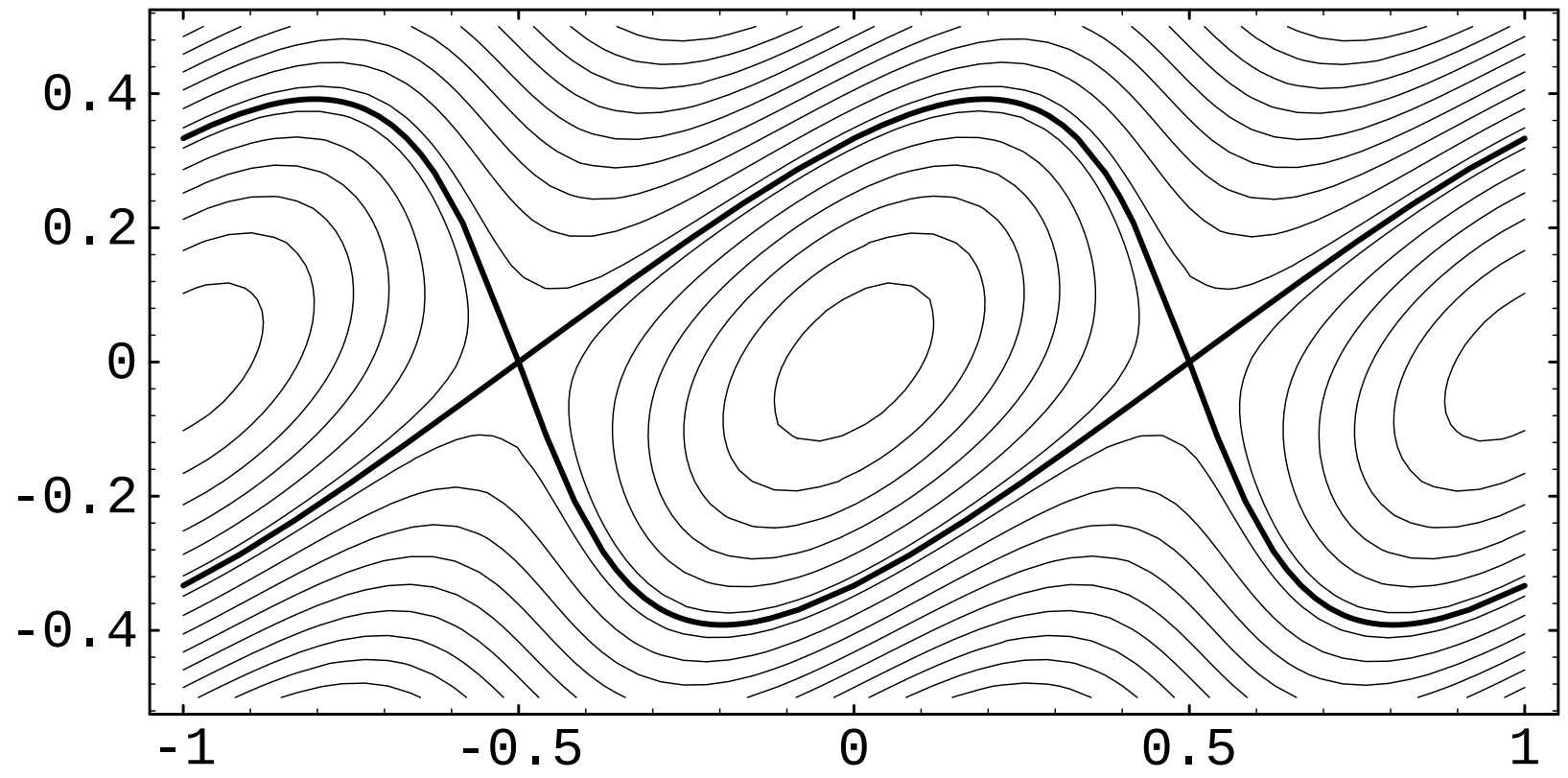

Figure 1: Some contours of $I_\delta$, for $\delta = 1/3$.

where the dilogarithm is defined by

$$\mathrm{dilog}(x) \equiv \int_1^x \frac{\log(z)}{1-z} dz \ .$$

The Suris map, $f_\delta$, i.e., the standard map with this potential, is integrable with integral

$$I_\delta(\theta, r) = \cos \pi r + \delta \cos \pi(2\theta - r) \ ,$$

i.e. $I_\delta \circ f_\delta = I_\delta$. Contours of $I_\delta$ are shown in figure 1. The map is integrable for any $\delta$, however, we will consider the case $0 < \delta < 1$, as the topology of the saddle connections changes at $\delta = 1$. For $0 < \delta < 1$ the map $f_\delta$ has hyperbolic fixed points at $z_a = (-\frac{1}{2}, 0)$ and $z_b = (\frac{1}{2}, 0)$, that are connected by two saddle connections, forming the upper and lower separatrices of the fixed point resonance.

Since the Suris map is integrable and has twist, Birkhoff's theorem [17] implies that the saddle connection between the equivalent points $z_a$ and $z_b$, which is a rotational invariant circle, is necessarily the graph of a function $r = \chi(\theta)$, for $\theta \in \bar{\mathcal{U}} = [-\frac{1}{2}, \frac{1}{2}]$. The dynamics on the saddle connection from (1), gives the circle diffeomorphism $\theta \to h(\theta) = \theta + \chi(\theta) + V'(\theta)$.

We will show that, under perturbation, the stable and unstable manifolds of the perturbed twist map intersect transversally. We use a modification of Melnikov's method [18] for twist maps [7, 8, 9].

We perturb the system by adding to the original potential any $C^2$ periodic function $P(\theta)$. We assume for simplicity that $P$ also satisfies $P(-\frac{1}{2}) = P(\frac{1}{2}) = 0$ and $P'(-\frac{1}{2}) = P'(\frac{1}{2}) = 0$. Then for small enough $\epsilon$, $S_\delta + \epsilon P$ is a twist generating function, and the corresponding map $f_{\delta,\epsilon}$ has hyperbolic fixed points at $z_a$ and $z_b$. As was shown in [8], the Melnikov series, given by

$$L(\theta) = \sum_{t=-\infty}^{\infty} P(h^t(\theta), h^{t+1}(\theta))$$

converges absolutely and uniformly to a $C^2$ function on compact subsets of $\mathcal{U}$. Furthermore, if the function $L$ has a nondegenerate critical point in $\mathcal{U}$, then the manifolds $W^u(z_a, f_{\delta,\epsilon})$ and $W^s(z_b, f_{\delta,\epsilon})$ intersect transversally for $\epsilon$ small enough (the same conclusion is valid for $W^s(z_a, f_{\delta,\epsilon})$ and $W^u(z_b, f_{\delta,\epsilon})$).

Thus the function $L$ provides sufficient conditions for the transversal destruction of the saddle connection, just as in the classic applications of the Melnikov integral. In section 2 we formulate a slightly stronger version of the previous result (theorem 2).

Below we consider the case $P(\theta) = \cos^2 \pi\theta$. We show that for $\epsilon > 0$ and small enough there are two distinct transversal heteroclinic orbits from $z_a$ to $z_b$. These delineate a turnstile with two lobes that define the areas that are transported into and out of the fixed point resonance upon each iteration of the map. The area of the lobe gives a coordinate independent measure of the separation of the perturbed stable and unstable manifolds as well as the flux from one region to another.

The exact, order $\epsilon$ approximation for the lobe area is obtained in section 4 and compared with numerical results in section 7. We summarize our results, as

**Theorem 1 (Main Theorem)** *Let*

$$V_\delta(\theta) = -\frac{2}{\pi} \int_0^\theta \arctan\left(\frac{\delta \sin(2\pi t)}{1 + \delta \cos(2\pi t)}\right) dt \tag{7}$$

*and* $S_\delta = \frac{1}{2}(\theta' - \theta)^2 + V_\delta(\theta)$ *be the generating function of the integrable Suris map* $f_\delta$. *Let*

$$S_{\delta,\epsilon} = S_\delta + \epsilon \cos^2 \pi\theta$$

*be a perturbed generating function with corresponding twist map* $f_{\delta,\epsilon}$. *Then for all* $0 < \delta < 1$ *and* $\epsilon$ *small enough there are two orbits heteroclinic from* $z_a = (-\frac{1}{2}, 0)$ *to* $z_b = (\frac{1}{2}, 0)$. *The lobe defined by these two orbits has an area* $A$ *given by*

$$A(\delta, \epsilon) = \epsilon\Gamma(\nu) + O(\epsilon^2) \tag{8}$$

*where*

$$\Gamma(\nu) \equiv 1 + 8\sum_{k=1}^{\infty} \frac{(-1)^k k\nu^k}{1+\nu^k}\ , \tag{9}$$

*and* $\nu \equiv (1-\sqrt{\delta})/(1+\sqrt{\delta})$.

Note that $\nu$ is a natural parameter for the stable and unstable manifolds, and hence for the area of the lobe, since the multipliers of the hyperbolic fixed point of the Suris map are $\nu$ and $1/\nu$.

The series $\Gamma(\nu)$ is rather intriguing. It is an analytic function of $\nu$ on the interval $\{0 < \nu < 1\}$, and approaches zero rapidly as $\nu$ increases. In fact, we will see in section 5 that $\Gamma(\nu)$ is strictly positive, but approaches zero exponentially fast:

$$\Gamma(\nu) \sim \left(\frac{4\pi}{\log(1/\nu)}\right)^2 \exp\left(\frac{-\pi^2}{\log(1/\nu)}\right)$$

as $\nu \to 1^-$. This, of course, corresponds to $\delta \to 0^+$. The analysis uses a remarkable relation between the series for $\Gamma$ and elliptic functions, and some formulas found in one of the notebooks of Ramanujan (cf. [2]).

## 2 Melnikov formula for twist maps

In this section we review derivation of the Melnikov formula for twist maps [8]. We begin with a $C^2$ Lagrangian generating function $S(\theta,\theta')$ that satisfies the twist condition (4) and has zero net flux, (5). It gives a map of the annulus implicitly through (3). Alternatively, the generating function defines the Lagrangian equations of motion through the action $W$ defined on a sequence $[\theta]=\{\theta^i,\theta^{i+1},\ldots,\theta^j\}$, by

$$W[\theta]=\sum_{t=i}^{j} S(\theta^t,\theta^{t+1})$$

It is easy to see that an orbit of the map that begins at $\theta^i$ and ends at $\theta^j$ corresponds to a critical point of the action $W$ under variation with respect to the interior points, and a periodic orbit is a critical point with respect to all points subject to the constraint that $\theta^j=\theta^i$. The corresponding momenta are then defined through (3) as $r^t=-\partial_1 S(\theta^t,\theta^{t+1})$. Thus, for example, a point $(\theta_a,r_a)$ a fixed point of the map if and only if $\theta_a$ is a critical point of $S(\theta,\theta)$, and $r_a$ is defined through (3).

We now describe how a modification of the *Melnikov method* can be applied to predict the transversal intersection between the stable and unstable manifolds of two different periodic or fixed points of a twist map. This method is based on the variational approach of Aubry [1] and Mather [15] and can be applied in general to any twist map that has a saddle connection, in particular to any integrable twist map. We begin with a map $f_0$ generated by $S_0$. Suppose that $f_0$ has two hyperbolic fixed points $z_a=(\theta_a,r_a)$ and $z_b=(\theta_b,r_b)$, and there is a saddle connection defined by the graph of a function $\chi(\theta)$ on the interval $\mathcal{U}=(\theta_a,\theta_b)$ between these points. A diffeomorphism $h:\bar{\mathcal{U}}\to\bar{\mathcal{U}}$ is induced by the restriction of the map to the saddle connection:

$$f_0(\theta,\chi(\theta))=(h(\theta),\chi(h(\theta)))$$

Let $P$ be a $C^2$ function with zero net flux, (5). Then the function

$$S_\epsilon(\theta,\theta')=S_0(\theta,\theta')+\epsilon P(\theta,\theta')$$

generates a twist map $f_\epsilon$ for small enough $\epsilon$. Since hyperbolic points are nondegenerate critical points of the action [10], the perturbed map will have nearby hyperbolic fixed points for small enough $\epsilon$. A simple case occurs when $\theta_a$ is a critical point of $P(\theta,\theta)$ as well as of $S_0(\theta,\theta)$ since it is then a critical point of $S_\epsilon(\theta,\theta)$ as well. Thus the fixed points of $S_\epsilon$ will have unchanged configurations, but their momenta will be modified according to (3).

It is well known that there are useful relations between the action of orbits and areas of regions for twist maps [14]. We will use one such relation to obtain the Melnikov formula: a relation, between the graph $\chi$ and the action of orbits on the stable and unstable manifolds of a hyperbolic fixed point [23]. Let $\theta^0$ be a point on the unstable manifold of a fixed point $z_a$ that is close enough to $z_a$ so that the segment of $W^u(z_a)$ containing $\theta^0$ is given by a graph $\chi(\theta)$. Let $\theta^t, t \le 0$, be the preorbit of this point. Defining the backward action difference as

$$\Delta W^b(\theta^0) = \sum_{t=-\infty}^{-1} \Big[ S(\theta^t, \theta^{t+1}) - S(\theta_a, \theta_a) \Big] \ ,$$

then the unstable manifold is defined by the graph of the function

$$\chi^u(\theta) = \frac{\partial \Delta W^b}{\partial \theta} \ .$$

A corresponding formula for the forward action of an orbit on the stable manifold yields a formula for the graph of an initial segment of the stable manifold, $\chi^s$, of $z_b$:

$$\begin{aligned} \Delta W^f(\theta^0) &= -\textstyle\sum_{t=0}^{\infty} \left[ S(\theta^t, \theta^{t+1}) - S(\theta_b, \theta_b) \right] \ , \\ \chi^s(\theta) &= \tfrac{\partial \Delta W^f}{\partial \theta} \ . \end{aligned} \qquad (10)$$

The difference between these two actions leads to the Melnikov-like formula for the transversal intersection of these manifolds. We summarize the results as a theorem:

**Theorem 2** *Let $S_0$ be the generating function for a twist map $f_0$ that has two hyperbolic fixed points $z_a$ and $z_b$ with a saddle connection given by $r = \chi(\theta)$ for $\theta \in \mathcal{U} = (\theta_a, \theta_b)$. Let $f_0$ induce a diffeomorphism $h(\theta)$ on the connection. Then if $S_\epsilon = S_0 + \epsilon P$ generates the twist map $f_\epsilon$, such that the perturbation $P$ has the following properties*

**a)** $P(\theta_a, \theta_a) = P(\theta_b, \theta_b) = 0$

**b)** $\frac{d}{d\theta}\Big|_{\theta=\theta_a} P(\theta,\theta) = \frac{d}{d\theta}\Big|_{\theta=\theta_b} P(\theta,\theta) = 0$

*Then for $\epsilon > 0$ small enough:*

**i)** *the perturbed map has two hyperbolic fixed points continued from $z_a$ and $z_b$;*

**ii)** *the Melnikov series*

$$L(\theta) = \sum_{t=-\infty}^{\infty} P(h^t(\theta), h^{t+1}(\theta)) \tag{11}$$

*converges absolutely and uniformly to a $C^2$ function on $\mathcal{U}$; and*

**iii)** *if $L$ has a nondegenerate critical point on $\mathcal{U}$, then the unstable and stable manifolds of the two fixed points intersect transversally.*

# 3 Transport

The theory of transport in two dimensional maps is based on a partition of phase space into regions between which transport is restricted by partial barriers of some sort. One of the simplest such partitions is to define a resonance zone associated with a saddle fixed point and its homoclinic tangle (or the heteroclinic tangle between two such points) [14, 20, 4]. We first recall a few definitions.

An *initial segment* of a (un)stable manifold of a saddle $z$ is defined as a segment starting at $z$ and continuing to some endpoint, say $p$. A *resonance zone* for $z$ is a compact region bounded by alternating initial segments of stable and unstable manifolds that intersect only at their endpoints. Each such intersection, apart from that at $z$ defines a *principal homoclinic point* $p$. Even with the choice of a particular homoclinic orbit there are an infinity of choices for a resonance zone, since each point $f^t(p)$ on the orbit of $p$ is a principal homoclinic point; however, for an area preserving map these regions all have the same area. To be physically meaningful, the point $p$ should be chosen so that the resonance zone is not too distorted; often symmetry can be used to give an appropriate choice.

We will see that, for a perturbation of the Suris map, there are two saddles $z_a$ and $z_b$ that have heteroclinic orbits going both from $z_a \to z_b$ and from $z_b \to z_a$, see figure 2. The resonance zone is made up from four initial segments, first a segment of $W^u(z_a)$ to $p$ that joins a segment of $W^s(z_b)$ at $p$ forming the upper part of the boundary, and then a segment of $W^u(z_b)$ to a lower homoclinic point, where it intersects a segment of $W^s(z_a)$ forming the lower piece of the boundary.

Escape from a resonance will be slow if the unstable and stable manifolds are nearly coincident. Recall that the *exit set* is the set that leaves the resonance upon one iteration, and the *incoming set* is the set that enters the resonance upon one iteration. These are easily obtained by taking the preimage of the boundary–all points on the unstable segments shrink towards their respective saddles, and the points on the stable segments lengthen. The preimage of a resonance is also a resonance, but the principal homoclinic points $p$ switch to $f^{-1}(p)$. For a principal point $p$ the segment of $W^s$ between $p$ and $f^{-1}(p)$ together with that of $W^u$ connecting these points bound a *turnstile* that is the union of the exit set and the incoming set.

For an area preserving map, the exit and incoming sets must have equal areas. In the simplest case there is one additional principal homoclinic point on the segment of the stable and unstable manifolds between $f^{-1}(p)$ and $p$, call it $q$. In this case the turnstile contains two *lobes*, the exit lobe and the incoming lobe. Specializing to the case of a twist map with generating function $S$, and denote the two principal homoclinic points by $p = (r_p, \theta_p)$ and $q = (r_q, \theta_q)$. In this case, the areas of the exit and entrance lobes are determined by the difference in action between the orbits of $p$ and of $q$ [13]

$$A = \sum_{t=-\infty}^{\infty} S(\theta_q^t, \theta_q^{t+1}) - S(\theta_p^t, \theta_p^{t+1}) \ . \tag{12}$$

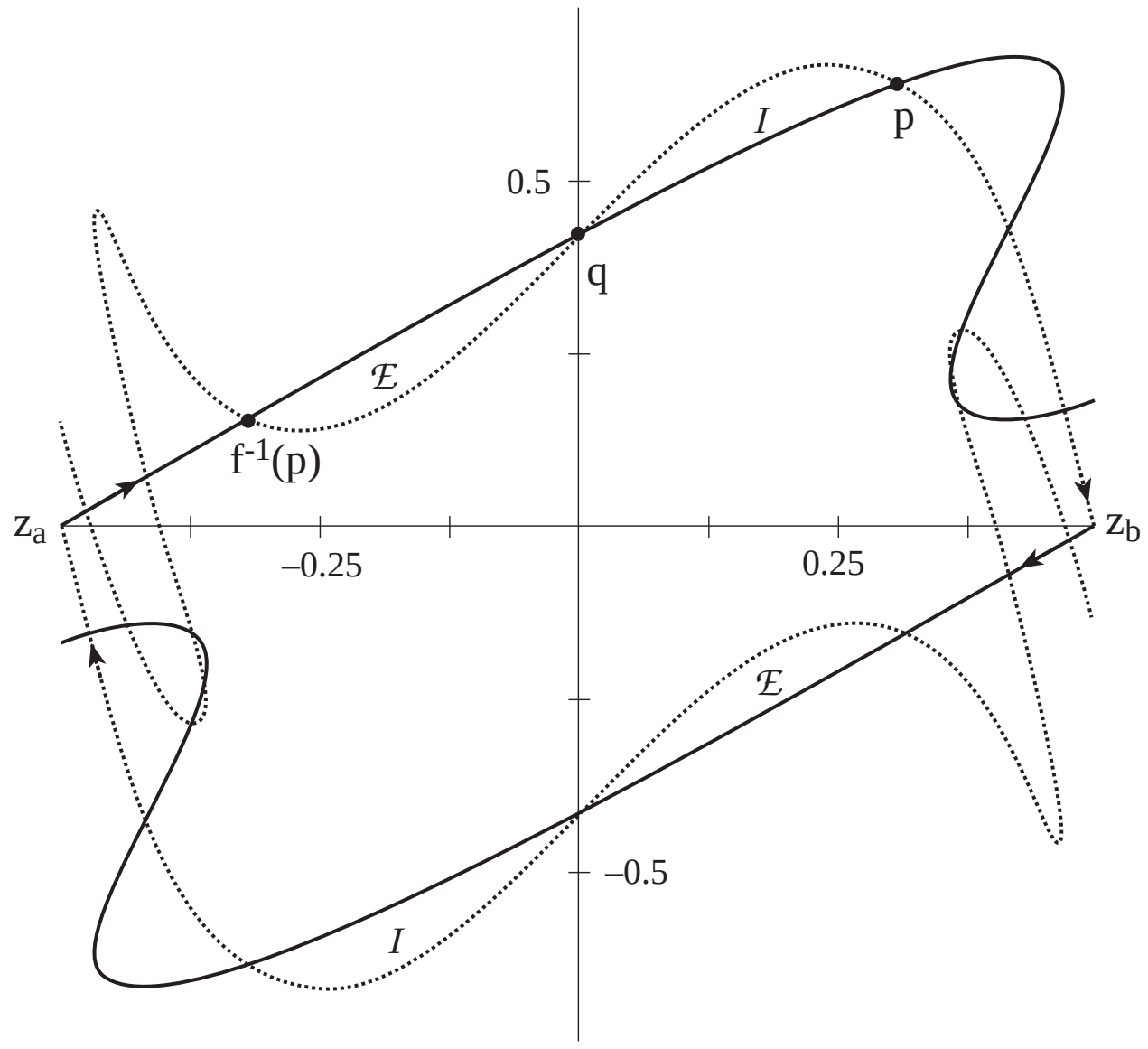

Figure 2: Resonance for the perturbed Suris map with $\delta = 0.5$ and $\epsilon = 0.05$. A pair of principal homoclinic points are labeled $p$ and $q$, and the exit and incoming sets are labeled $\mathcal{E}$ and the $\mathcal{I}$.

The area $A$ is the signed area below the segment of $W^s$ between $q$ and $p$ and above that of $W^u$.

Now, the action is stationary on an orbit, and when there is a saddle connection, the actions of $p$ and $q$ are equal. Thus it is easy to see that the area of the exit set for a perturbation of a map with a saddle connection is determined, to lowest order, by the difference between two critical values of the Melnikov function $L$ [11]. We summarize this result as a theorem.

**Theorem 3** *Let $S_0$, $P$, $S_\epsilon$, $f_0$, $f_\epsilon$, $\mathcal{U}$, $h$ and $L$ as in section 2. Assume that $\theta_p$ and $\theta_q$ are two points in $\mathcal{U}$ such that*

**a)** $L'(\theta_q) = L'(\theta_p) = 0$

**b)** $L''(\theta_p) > 0$, $L''(\theta_q) < 0$

**c)** $L'(\theta) \neq 0$ *for* $\theta_q < \theta < \theta_p$

*Then the heteroclinic points $p$ and $q$ continue heteroclinic points of $f_\epsilon$, and the stable and unstable manifolds of $f_\epsilon$ enclose a lobe with area*

$$A(\epsilon) = \epsilon(L(\theta_q) - L(\theta_p)) + O(\epsilon^2)$$

# 4 The Suris Map

In this section we give formulas for the dynamics on the homoclinic connection of the fixed point. Recall that the Suris map we consider is generated by the function (6) with the potential (7)

Since $V$ has minima at $\theta = \frac{n}{2}$, the Suris map has hyperbolic fixed points at $(\theta, r) = (-\frac{n}{2}, 0)$. In fact, there are saddle connections between neighboring points and their exact formulas can be computed. Setting $\mathcal{U} = \{-\frac{1}{2} < \theta < \frac{1}{2}\}$, we find that the diffeomorphism on the saddle connection, $h_\nu : \bar{\mathcal{U}} \to \bar{\mathcal{U}}$ is

$$h_\nu(\theta) = \frac{2}{\pi} \arctan\left( \frac{(\nu+1)\tan(\frac{\pi}{2}\theta) + (\nu-1)}{(\nu-1)\tan(\frac{\pi}{2}\theta) + (\nu+1)} \right) , \tag{13}$$

where $\nu = (1-\sqrt{\delta})/(1+\sqrt{\delta})$. Then we have the following.

**Lemma 4** *Let $h_\nu$ be the diffeomorphism of $\bar{\mathcal{U}}$ given by (13). Then $h_\nu$ satisfies*

**a)** $h_\nu^t = h_{\nu^t}$, *for all* $t \in \mathbb{Z}$.

**b)** $h_\nu(-\theta) = -h_{\nu^{-1}}(\theta)$.

**c)** $-\frac{1}{2}$ *is an unstable fixed point and* $\frac{1}{2}$ *is a stable fixed point of* $h_\nu$.

**d)** $V_\delta'(\theta) = h_\nu(\theta) - 2\theta + h_\nu^{-1}(\theta)$, *for* $\theta \in \bar{\mathcal{U}}$, *for* $V_\delta$ *given by (7).*

**Proof.** A direct computation proves a) and b).

Direct computation for c) shows that the only fixed points of $h_\nu$ occur when $\tan(\frac{\pi}{2}\theta) = \pm 1$, or $\theta = \pm\frac{1}{2}$. Furthermore, $h_\nu'(-\frac{1}{2}) = \frac{1}{\nu}$ and $h_\nu'(\frac{1}{2}) = \nu$. Since $\nu \in (0,1)$ this implies that the former is unstable and the latter is stable.

For d), we will first show that when $\theta \in \bar{\mathcal{U}}$,

$$\theta - h_\nu(\theta) = \frac{2}{\pi} \arctan\left( \frac{\sqrt{\delta}\cos\pi\theta}{1-\sqrt{\delta}\sin\pi\theta} \right) , \tag{14}$$

and

$$\theta - h_\nu^{-1}(\theta) = \frac{2}{\pi} \arctan\left( \frac{-\sqrt{\delta}\cos\pi\theta}{1+\sqrt{\delta}\sin\pi\theta} \right) . \tag{15}$$

In order to prove these formulas, we notice first that if $\theta \in \mathcal{U}$, then $h(\theta) \in \mathcal{U}$. Therefore $-1 < \theta - h_\nu(\theta) < 1$. On the other hand,

$$\tan\left(\frac{\pi}{2}[\theta - h_\nu(\theta)]\right) = \frac{\sqrt{\delta}\cos\pi\theta}{1-\sqrt{\delta}\sin\pi\theta} \ , \tag{16}$$

This implies (14). It is easy to see that if $\theta \in \mathcal{U}$, then (16) is positive and therefore $0 < \theta - h_\nu(\theta) < 1$. The substitution $-\theta \mapsto \theta$ in equation (14), implies that $0 < -\theta + h_{\nu^{-1}}(\theta) < 1$, and equation (15). We conclude that for all $\theta \in \mathcal{U}$, $-1 < 2\theta - h_{\nu^{-1}}(\theta) - h_\nu(\theta) < 1$. To finish the proof we work on the tangent of the second difference

$$\begin{aligned}
&\tan\left(\frac{\pi}{2}[2\theta - h_{\nu^{-1}}(\theta) - h_\nu(\theta)]\right) = \\
&\tan\left(\tan^{-1}\left(\frac{-\sqrt{\delta}\cos\pi\theta}{1+\sqrt{\delta}\sin\pi\theta}\right) + \tan^{-1}\left(\frac{\sqrt{\delta}\cos\pi\theta}{1-\sqrt{\delta}\sin\pi\theta}\right)\right) \\
= \ &\frac{\sqrt{\delta}\cos\pi\theta(2\sqrt{\delta}\sin\pi\theta)}{1-\delta\sin^2\pi\theta+\delta\cos^2\pi\theta} \\
= \ &\frac{\delta\sin(2\pi\theta)}{1+\delta\cos(2\pi\theta)}
\end{aligned}$$

Since this last expression is $-\tan(\frac{\pi}{2}V_\delta'(\theta))$, this completes the proof. □

## 4.1 Saddle Connections

With the help of lemma 4, we can give a description of the intersection of the saddle connection between $(-\frac{1}{2}, 0)$ and $(\frac{1}{2}, 0)$. The saddle connection is shown in figure 3.

**Proposition 5** *Let $f_\delta$ be the twist map generated by $S_\delta$. Then $(-\frac{1}{2}, 0)$ and $(\frac{1}{2}, 0)$ are hyperbolic fixed points for $f_\delta$ and there are two saddle connections between them given by the graphs*

$$r = \chi^+(\theta) = \theta - h_{\nu^{-1}}(\theta) \ ,$$

*and*

$$r = \chi^-(\theta) = \theta - h(\theta) \ .$$

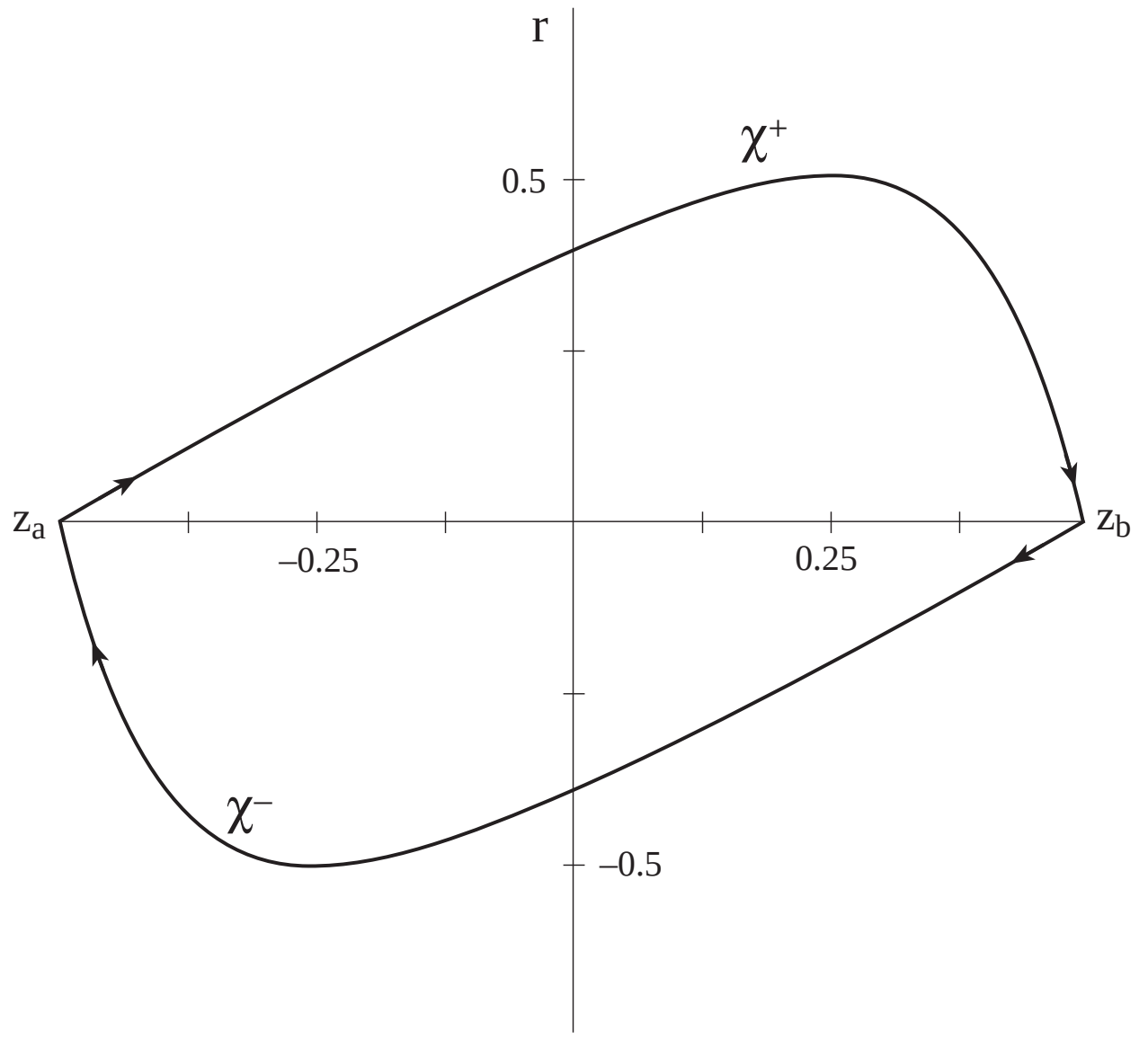

Figure 3: Saddle connections for the Suris map with $\delta = 1/2$.

**Proof.** Using Lemma 4, the map can be written in the form

$$f(\theta, r) = (r + h_\nu(\theta) - \theta + h_{\nu^{-1}}(\theta), r + h_\nu(\theta) - 2\theta + h_{\nu^{-1}}(\theta))$$

Therefore, if $r = \chi^+(\theta) = \theta - h_{\nu^{-1}}(\theta)$ then

$$f(\theta, \chi^+(\theta)) = (h_\nu(\theta), h_\nu(\theta) - \theta) = (h_\nu(\theta), \chi^+(h_\nu(\theta)) \ .$$

Now since $\theta = -\frac{1}{2}$ is an unstable fixed point for $h_\nu$, this graph gives the right going saddle connection. In the same way setting $r = \chi^-(\theta) = \theta - h_\nu(\theta)$ gives

$$f(\theta, \chi^-(\theta)) = (h_{\nu^{-1}}(\theta), h_{\nu^{-1}}(\theta) - \theta) = (h_{\nu^{-1}}(\theta), \chi^-(h_{\nu^{-1}}(\theta)) \ .$$

Since the map conjugates to $h_{\nu^{-1}}$ on this graph, this is clearly the left going saddle connection. □

## 4.2 Proof of the Main Theorem

In this section we sketch the proof of the main theorem. The analysis of the infinite series for the Melnikov function relies on some formulae that can be found in one of the notebooks of Ramanujan (cf. [2]), though we are not sure if he is the original author of the formulas or the first to publish them. However, the formulas are remarkable and all that we do in this section and the next is a consequence of them. First we recall the theorem:

**Theorem 6 (Main Theorem)** *Let $V_\delta$ be given by (7) and $S_\delta$ be the Suris generating function (6). Let*

$$S_{\delta,\epsilon} = S_\delta + \epsilon \cos^2 \pi\theta$$

*be a perturbed generating function with corresponding twist map $f_{\delta,\epsilon}$.*

*Then for $0 < \delta < 1$ there are at least two distinct transversal heteroclinic orbits connecting $(-1/2, 0)$ and $(1/2, 0)$. Furthermore, the stable and unstable manifolds of these fixed points enclose a lobe with area $A(\delta,\epsilon)$ given by*

$$A(\delta,\epsilon) = \epsilon\left(1 + 8\sum_{k=1}^{\infty} \frac{(-1)^k k\nu^k}{1+\nu^k}\right) + O(\epsilon^2) \ ,$$

*where $\nu = (1-\sqrt{\delta})/(1+\sqrt{\delta})$.*

**Proof.** Let $P(\theta) = \cos^2 \pi\theta$. It is clear that $P$ satisfies the conditions of Theorem 2. Let

$$L(\theta) = \sum_{t=-\infty}^{\infty} P(h^t(\theta)) \ .$$

According to theorem 2 a sufficient condition for transversal intersection of the perturbed manifolds is that $L$ has a nondegenerate critical point on the interval $\mathcal{U} = \{-\frac{1}{2} < \theta < \frac{1}{2}\}$. The graph of $L$ over $\mathcal{U}$ is shown in figure 4.

To proceed, we define $\bar{L}(z)$ by $L(\theta) = \bar{L}(\tan(\frac{\pi}{2}\theta))$; therefore $\bar{L}(z) = \sum_{t=-\infty}^{\infty} \alpha_t(z)$ where

$$\alpha_t(z) = \frac{4\nu^{2t}(1-z^2)^2}{((1-z)^2 + \nu^{2t}(1+z)^2)^2} \ .$$

This implies that $\theta$ is a non degenerate critical point of $L$ whenever $\tan(\frac{\pi}{2}\theta)$ is a nondegenerate critical point of $\bar{L}$. We are going to show that $\theta_q = 0$ is a local maximum of $L$ and $\theta_p = \frac{2}{\pi}\arctan((1-\sqrt{\nu})/(1+\sqrt{\nu}))$ is a local minimum.

First, we write $\bar{L}$ as

$$\bar{L}(z) = \alpha_0(z) + \sum_{t=1}^{\infty}\{\alpha_t(z) + \alpha_{-t}(z)\} \ . \tag{17}$$

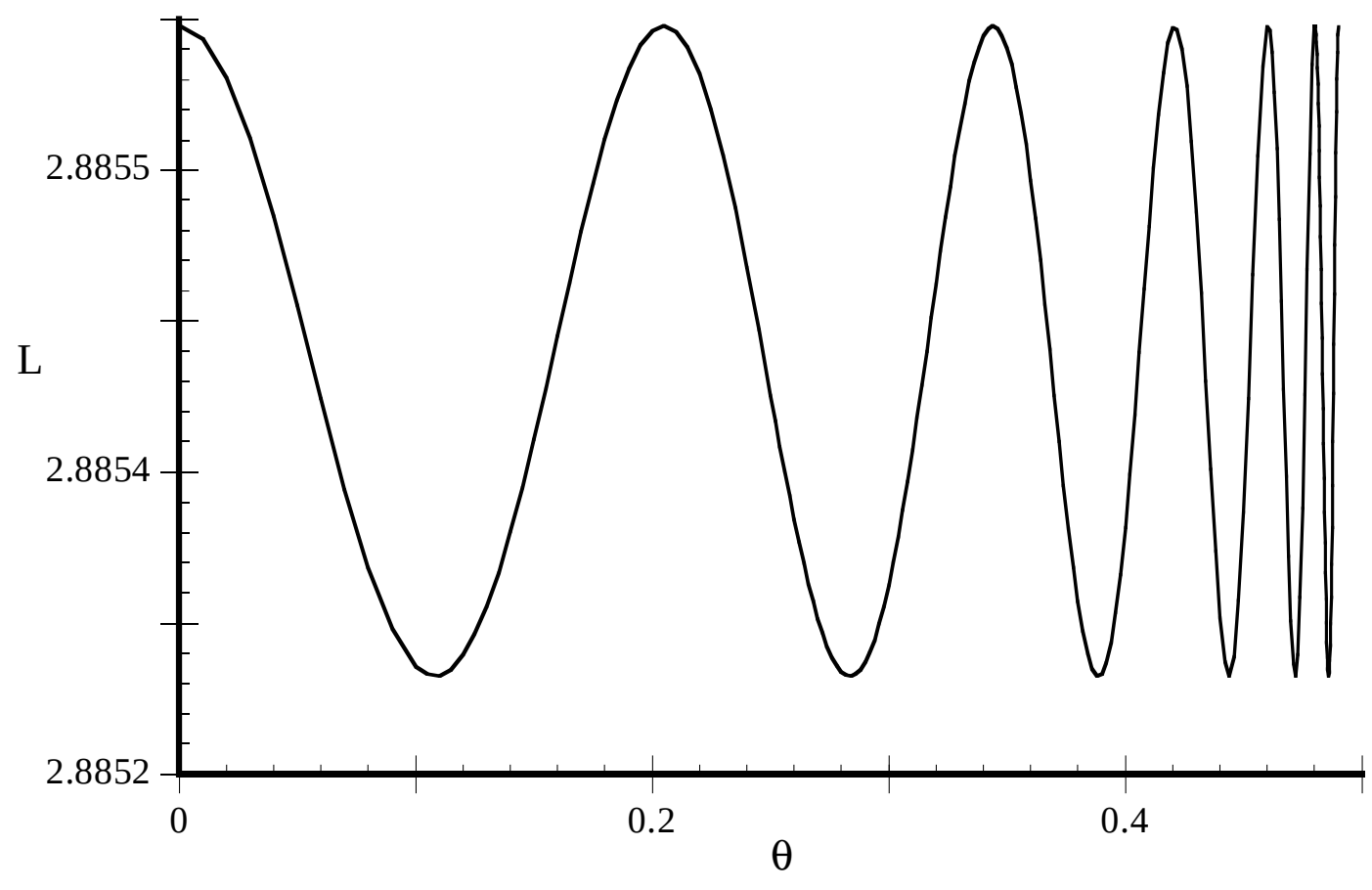

Figure 4: Melnikov function, $L$ for $\delta = 1/2$.

Now, since $\alpha_t(-z) = \alpha_{-t}(z)$, this implies that $\bar{L}$ is even, and so $\bar{L}'(0) = 0$. Differentiation gives

$$\bar{L}''(0) = -8 + \sum_{t=1}^{\infty} \frac{128\left(1 - 4\nu^{2t} + \nu^{4t}\right)\nu^{2t}}{\left(1 + \nu^{2t}\right)^4} \quad .$$

Let

$$\Gamma_0(\nu) = 1 - 16 \sum_{t=1}^{\infty} \frac{\left(1 - 4\nu^{t} + \nu^{2t}\right)\nu^{t}}{\left(1 + \nu^{t}\right)^4} \quad ,$$

then $\bar{L}''(0) = -8\Gamma_0(\nu^2)$. Expanding the denominator of $\Gamma_0$ and rearranging the sums gives

$$\Gamma_0(\nu) = 1 + 16 \sum_{k=0}^{\infty} \frac{(-1)^k k^3 \nu^k}{1 - \nu^k} \quad .$$

Now let $K(x)$ be the normalized complete elliptic integral of the first kind

$$K(x) = \frac{2}{\pi} \int_0^{\pi/2} \frac{d\tau}{\sqrt{1 - x \sin^2 \tau}} \quad . \tag{18}$$

and define an increasing diffeomorphism of the interval $(0, 1)$

$$H(x) = \exp\left(-\pi \frac{K(1-x)}{K(x)}\right) \quad . \tag{19}$$

Then, a Ramanujan formula [2, III.17, entry 14 (v)] implies that

$$\Gamma_0(H(x)^2) = K(x)^4(1-x) \ ,$$

and therefore $\bar{L}''(0) < 0$, for all $0 < \nu < 1$. This shows that $\theta_q = 0$ is a nondegenerate local maximum of $L$, for all $0 < \nu < 1$.

To find the local minimum, we rewrite $\bar{L}$ as

$$\bar{L}(z) = \sum_{t=1}^{\infty} \{\alpha_{t+1}(z) + \alpha_{-t}(z)\} \ . \tag{20}$$

Now we find that $\bar{L}'((1-\sqrt{\nu})/(1+\sqrt{\nu})) = 0$ and

$$\bar{L}''((1-\sqrt{\nu})/(1+\sqrt{\nu})) = \sum_{t=0}^{\infty} \frac{\left(1 - 4\nu^{1+2t} + \nu^{2+4t}\right)\nu^{1+2t}}{(1+\nu^{1+2t})^4} \ .$$

Notice that $\bar{L}''((1-\sqrt{\nu})/(1+\sqrt{\nu})) = \Gamma_0(\nu^2) - \Gamma_0(\nu)$. After the change of coordinates $\nu = H(x)$, we get from formulas (III.17.14(v)) and (III.17.14(ix)) in [2] that

$$\begin{aligned} \Gamma_0(H(x)^2) - \Gamma_0(H(x)) &= K(x)^4(1-x) - K(x)^4(1-x)^2 \\ &= K(x)^4(1-x)x \ , \end{aligned}$$

and therefore $\bar{L}''((1-\sqrt{\nu})/(1+\sqrt{\nu})) > 0$, for all $0 < \nu < 1$. This shows that $\theta_p$ is a nondegenerate local minimum of $L$, for all $0 < \nu < 1$.

Using Theorem 3, we conclude that a lobe of area

$$A(\delta, \epsilon) = \epsilon\,(L(\theta_q) - L(\theta_p)) + O(\epsilon^2)$$

is enclosed by the stable and unstable manifolds, where $\theta_p$ and $\theta_q$ are given above. Finally, we use (17,20) to obtain

$$L(\theta_q) = \bar{L}(0) = 1 + 8\sum_{t=1}^{\infty} \frac{\nu^{2t}}{(1+\nu^{2t})^2} \ ,$$

and

$$L(\theta_p) = \bar{L}((1-\sqrt{\nu})/(1+\sqrt{\nu})) = 8\sum_{t=0}^{\infty} \frac{\nu^{2t+1}}{(1+\nu^{2t+1})^2} \ .$$

Therefore,

$$\begin{aligned} L(\theta_q) - L(\theta_p) &= 1 + 8\sum_{t=1}^{\infty}\left(\sum_{k=1}^{\infty}(-1)^{k-1}k(\nu^{2t})^k\right) \\ &\quad -8\sum_{t=1}^{\infty}\left(\sum_{k=1}^{\infty}(-1)^{k-1}k(\nu^{2t-1})^k\right) \\ &= 1 + 8\sum_{k=1}^{\infty}\frac{(-1)^k k\nu^k}{1+\nu^k} = \Gamma(\nu) \ , \end{aligned}$$

where $\Gamma(\nu)$ was defined in (9). This concludes the proof. □

The explicit formula for the first order approximation of the area is compared with numerical computations in Sec. 7. We conclude with a corollary that will be useful in the next section.

**Corollary 7** *Let $A(\delta,\epsilon)$ be the area of the lobe that was described before. Let $K(x)$ and $H(x)$ as in the proof of the Main Theorem.*

*Let $G(x) = ((1 - H(x))/(1 + H(x))^2$. Then*

$$\lim_{\epsilon\to 0}\frac{A(G(x),\epsilon)}{\epsilon} = K(x)^2(1-x)$$

**Proof.** If $\delta = G(x)$ then $\nu = (1-\sqrt{\delta})/(1+\sqrt{\delta}) = H(x)$. Using [2, formula (III.17.14.(i))] we find that

$$\Gamma(H(x)) = K(x)^2(1-x). \tag{21}$$

□

# 5 Exponentially small behavior

In this section we investigate the asymptotics of the lobe area, (8), as $\delta \to 0$. We will see that the area is exponentially small in $\delta^{-1/2}$. The analysis in this section is also based on the Ramanujan formalae used in the previous section. We summarize with a lemma.

**Lemma 8** *Let $\Gamma(\nu)$ be defined by (9), then*

$$\Gamma(\nu) \sim \left(\frac{4\pi}{\log(1/\nu)}\right)^2 \exp\left(\frac{-\pi^2}{\log(1/\nu)}\right) \quad \textit{as } \nu \to 1^- \ . \tag{22}$$

**Proof.** We use (21) of corollary 7, based on rewriting $\Gamma$ in terms of the elliptic integral (18) and the diffeomorphism $H$, (19). Our elliptic function is normalized so that $K(0) = 1$; it has the asymptotic form [3, formula 112.01]

$$K(x) \sim \frac{1}{\pi} \log(\frac{16}{1-x}) \quad \text{as } \nu \to 1^- \ . \tag{23}$$

Thus $H$ has the limits

$$\lim_{x \to 0^+} H(x) = 0 \ , \quad \lim_{x \to 1^-} H(x) = 1 \ .$$

Further, (23) gives

$$H(x) \sim \frac{x}{16} \quad \text{as } x \to 0^+ \Longrightarrow H^{-1}(\nu) \sim 16\nu \quad \text{as } \nu \to 0^+ \ . \tag{24}$$

The definition of $H$ implies that $\log H(x) \log H(1-x) = \pi^2$ and therefore that

$$H(1-x) = \exp\left(\frac{\pi^2}{\log(\nu)}\right) = \exp\left(-\frac{\pi^2}{\log(1/\nu)}\right) \ ;$$

and so

$$1 - x = 1 - H^{-1}(\nu) = H^{-1}\left(\exp(-\pi^2/\log(1/\nu))\right) \ .$$

Combining this with (24) gives

$$1 - x \sim 16 \exp(-\pi^2/\log(1/\nu)) \text{ as } \nu \to 1^- \ . \tag{25}$$

Finally, using (25) and (23) yields

$$K(x) = K(H^{-1}(\nu)) \sim \frac{\pi}{\log(\frac{1}{\nu})} \quad \text{as } \nu \to 1^- \ . \tag{26}$$

Putting equations (25) and (26) into (21) gives the promised result. □

The asymptotic result for $\nu \to 1^-$ can be easily converted into one for $\delta \to 0^+$ using

$$\frac{1}{\log(1/\nu)} \sim \frac{1}{2\sqrt{\delta}}$$

Combining this with lemma 8 gives

**Corollary 9**

$$A(\delta, \epsilon) \sim \frac{\epsilon 4\pi^2}{\delta} \exp\left(\frac{-\pi^2}{2\sqrt{\delta}}\right) + O(\epsilon^2) \ , \tag{27}$$

*as* $\delta \to 0$.

## 6 Anti-Integrable Limit

While it doesn't fit in with the rest of our analysis, we present in this section a large $\epsilon$ expression for the lobe area. Our main purpose is to have an expression to compare with the numerical results in the next section. This expression is easy to obtain using the "anti-integrable limit" [12]. This limit is obtained by scaling the action by $\epsilon^{-1}$, to get

$$\hat{S} = S(\theta, \theta')/\epsilon = P(\theta) + \epsilon^{-1} S_\delta \ ,$$

and then setting $\epsilon^{-1} = 0$. The point is that when $\epsilon \gg 1 > \delta$ the points on the two heteroclinic orbits are all found in a neighborhood of the critical points of the potential $P = \cos^2(\pi\theta)$, i.e., at $\theta = m/2$. In the anti-integrable limit, an orbit consists solely of a sequence of configuration points sitting at these critical points. Thus the two heteroclinic orbits are given by the sequences

$$\{\theta_p^t\} = \{..., -\frac{1}{2}, -\frac{1}{2}, \frac{1}{2}, \frac{1}{2}, ...\} \ \ , \ \ \{\theta_q^t\} = \{..., -\frac{1}{2}, -\frac{1}{2}, 0, \frac{1}{2}, \frac{1}{2}, ...\}$$

The action is stationary on an orbit, thus to first order in the small parameter $\epsilon^{-1}$, the change in the orbit with $\epsilon$ can be ignored in the action difference to give

$$\frac{1}{\epsilon} A = \sum_{t=-\infty}^{\infty} \left( \hat{S}(\theta_q^t, \theta_q^{t+1}) - \hat{S}(\theta_p^t, \theta_p^{t+1}) \right) + O(\epsilon^{-2}) \ .$$

For the Suris map, this yields

$$A = \epsilon - \frac{1}{4} - \frac{1}{\pi^2} \left[ \mathrm{dilog}(1+\delta) - \mathrm{dilog}(1-\delta) \right] + O(\epsilon^{-1}) \qquad (28)$$

We compare this result with the numerical calculations in section 7.

# 7 Numerical comparison

In this section, we compare the theoretical results with numerical computations of the lobe area. The task is to find the actions of the two homoclinic points $p$ and $q$. For this task we use the symmetry of the Suris map. The action difference between the orbits of $q$ and $p$ gives the lobe area, (12).

A map is *reversible* if it is conjugate to its inverse by an involution $R$: $R^2 = I$ and $RfR = f^{-1}$. We call $R$ a reversor for $f$. Providing $V$ is an even function, a reversor for the standard map is given by

$$R: \ (\theta, r) \ \rightarrow \ (-\theta, r \ + V'(\theta)) \ . \qquad (29)$$

Note that if R is a reversor, then $f^t R$ is as well. Fixed sets, $\mathrm{Fix}(R) = \{z : Rz = z\}$ of a reversor are important because

**Lemma 10** *Let $z_a$ be a saddle fixed point of $f$. Suppose there is a point $q \in W^u(z_a)$ that is fixed under $R$, then the orbit of $q$ is heteroclinic from $z_a$ to $Rz_a$.*

**Proof.** Denote the orbit of $q$ by $q^t$, with $q = q^0$. By assumption $\lim_{t\to\infty} q^{-t} = z_a$. Using the reversor gives $q^t = f^t q^0 = f^t R q^0 = R f^{-t} q^0 = R q^{-t}$, thus $\lim_{t\to\infty} q^t = Rz_a$. □

Thus to find a heteroclinic orbit from $z_a$ to $Rz_a$ it is sufficient to search for points on the unstable manifold that intersect the fixed set of $R$. Such an orbit is called a *symmetric* heteroclinic orbit. Note that if $q^0$ is fixed under $R$, then $q^t \in \mathrm{Fix}(f^{2t}R)$ since

$$f^{2t} R q^t = f^t R f^{-t} q^t = f^t R q^0 = f^t q^0 = q^t \ .$$

This implies that the symmetric heteroclinic orbits divide naturally into pairs, one has a point on $\mathrm{Fix}(R)$, which we will call $q$ and the other on $\mathrm{Fix}(fR)$, called $p$. The fixed sets of the standard symmetry, (29), are

$$\mathrm{Fix}(R) = \{(\theta, r) : \theta = 0\} \ , \quad \mathrm{Fix}(fR) = \{(\theta, r) : r = 2\theta\} \ ,$$

see figure 5. To find these orbits numerically we move along $W^u(z_a)$ to the first points that intersect $\mathrm{Fix}(R)$ and $\mathrm{Fix}(fR)$, respectively. For example, to find $q \in \mathrm{Fix}R$, let $z^0 = (\theta^0, r^0) = z_a + \eta E^u$ where $E^u$ is the unstable eigenvector, and $\eta$ is a small parameter to be chosen below. Let $t_c + 1$ be the first time for which the iterate of $z^0$ is beyond $\mathrm{Fix}(R)$. Now choose a point $z(\theta)$ on the line from $z^0$ to $z^1 = fz^0$, parameterized by the initial angle $\theta$. By construction we are guaranteed that the function $Z(\theta) = \pi_1 f^{t_c}(z(\theta))$ has a zero for $\theta \in [\theta^0, \theta^1]$.

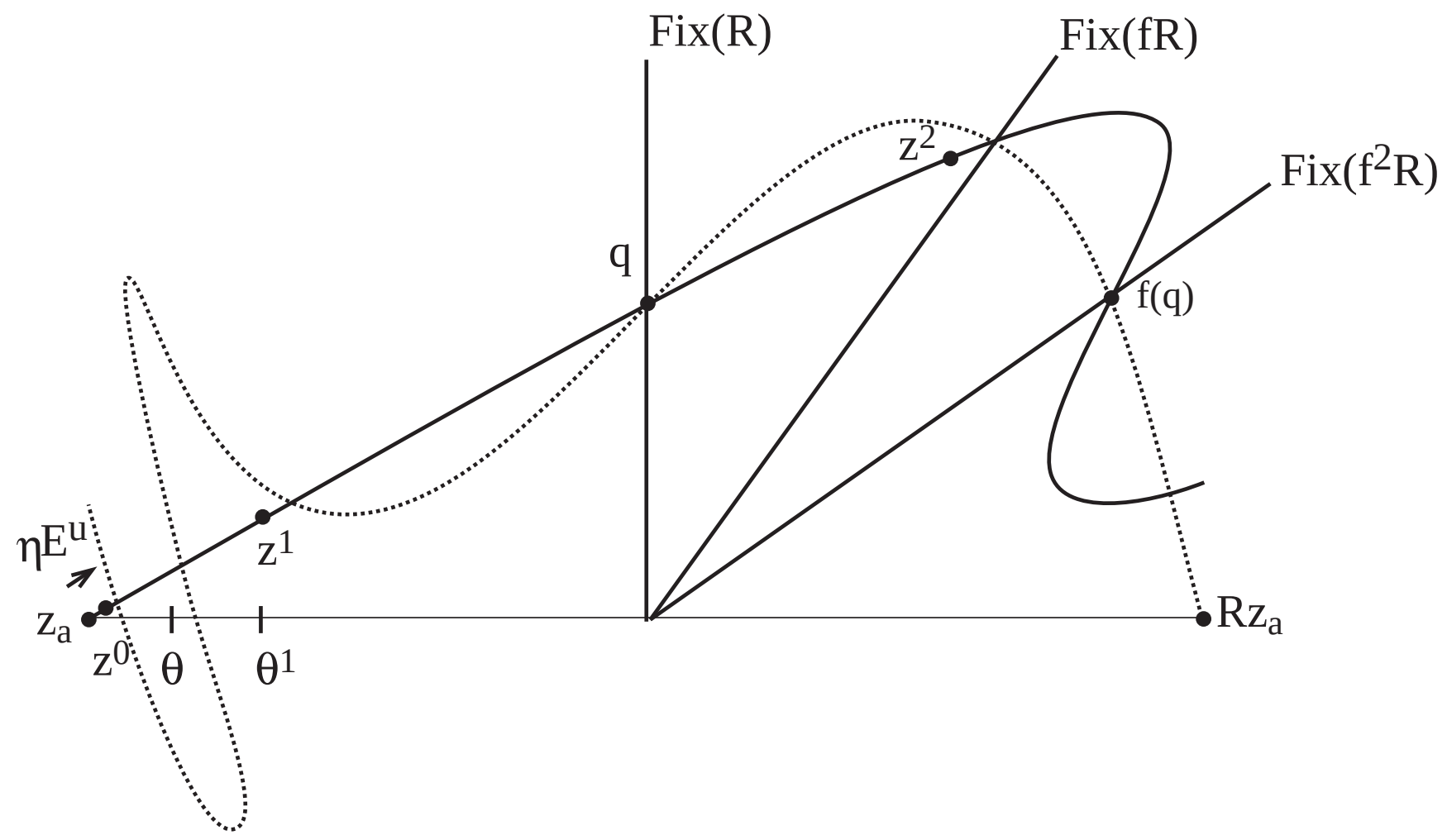

Figure 5: Finding some heteroclinic orbits using symmetry. Shown are the fixed sets $\mathrm{Fix}(f^t R)$, for $t = 0, 1, 2$ and a rough example of initial guess $z_0$, leading to $t_c = 1$.

We use a root finding method (Brent's method) to determine this zero to some precision, say $\rho$. The choice of precision influences the original value for $\eta$, as well as the number of iterates until a crossing. Assuming $W^u$ is smooth, the point $z^0$ will be $O(\eta^2)$ away from $W^u$. After $t_c$ iterates, however, this error will decrease by the factor $\lambda^{-t_c}$ where $\lambda$ is roughly the unstable multiplier of the fixed point. There is no sense in having this error smaller or larger than the precision of our root finder, so we set $\rho \sim \lambda^{-t_c}\eta^2$. On the other hand, since we start a distance $\eta$ from the fixed point, and wish to go a distance $O(1)$ to find the first crossing of the symmetry line, we have $\eta\lambda^{t_c} = O(1)$. Thus, it is appropriate to set

$$\eta \sim \rho^{1/3} \tag{30}$$

To find the second homoclinic point, $p \in \mathrm{Fix}(fR)$, we repeat the above analysis, using crossing of $\mathrm{Fix}(fR)$ to determine $t_c$, etc. The lobe area is given by the difference in action between these two orbits, from (12).

For our computations, using IEEE double precision arithmetic, we set

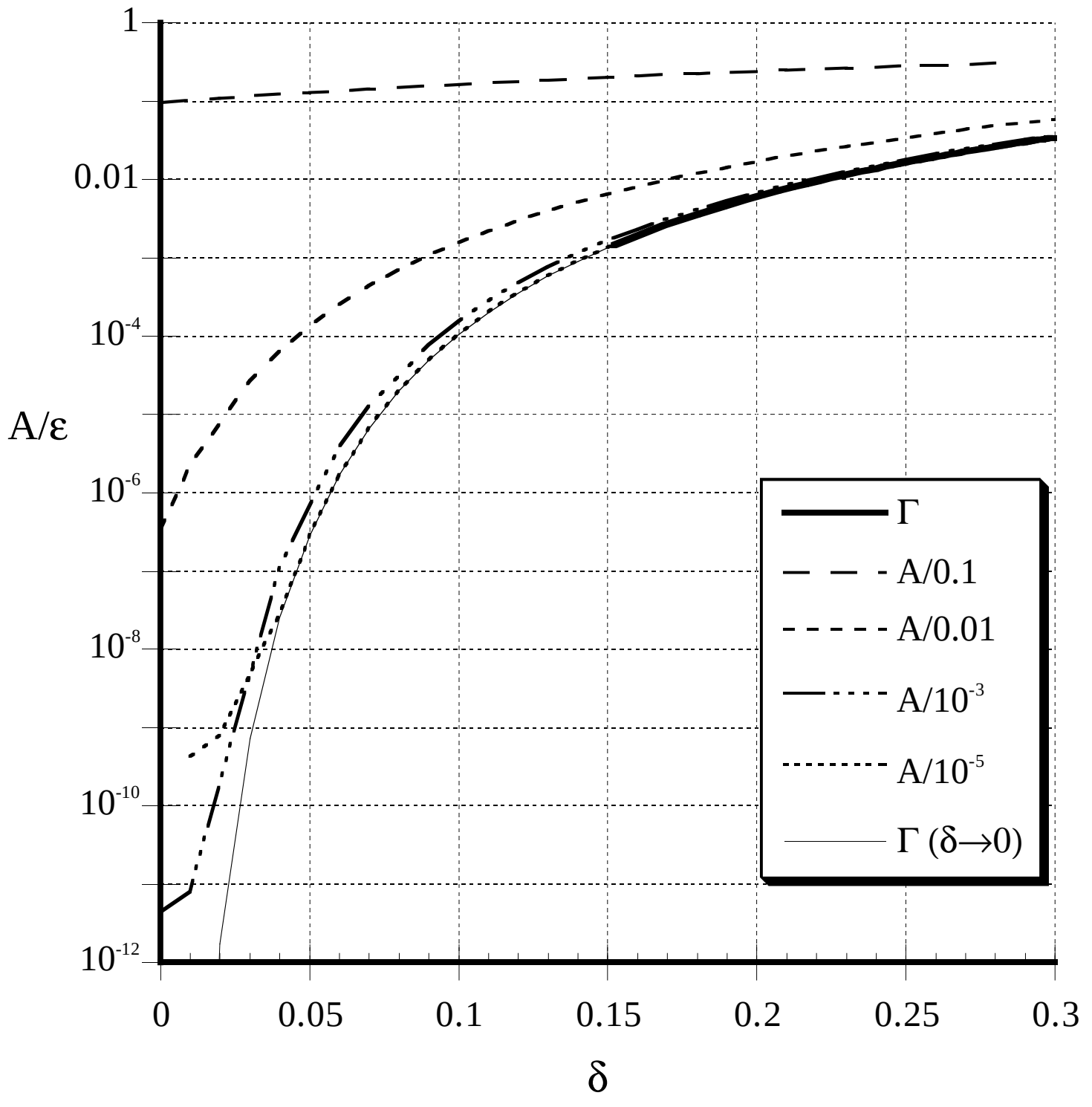

Figure 6: Log-Linear plot of $A/\epsilon$ as a function of $\delta$ for various values of $\epsilon$ compared with the theoretical expression (8) (thick line) and the small $\delta$ expression (22) (thin line)

$\rho = 10^{-19}$. These computations give apparently accurate results providing $A >> 10^{-14}$.

In figure 6 we show a comparison of the result of (8) with the numerical results on a log scale. The analytical result agrees well with the numerical results when $\epsilon = 10^{-5}$. We show the same data on a linear scale in figure 7. Remarkably, the asymptotic formula (22) agrees with the $\epsilon = 10^{-5}$ computation within 1% up to $\delta = 0.8$. Note, however, that (27) provides a poor comparison with the numerical results since even when $\delta = 0.1$, $\log(1/\nu)$ differs from $2\sqrt{\delta}$ by almost 4%. In figure 7, the anti-integrable results evaluated at $\epsilon = 1$ are also shown. We are unable to obtain numerical results for such a large $\epsilon$, as the multiplier of the fixed points is too large.

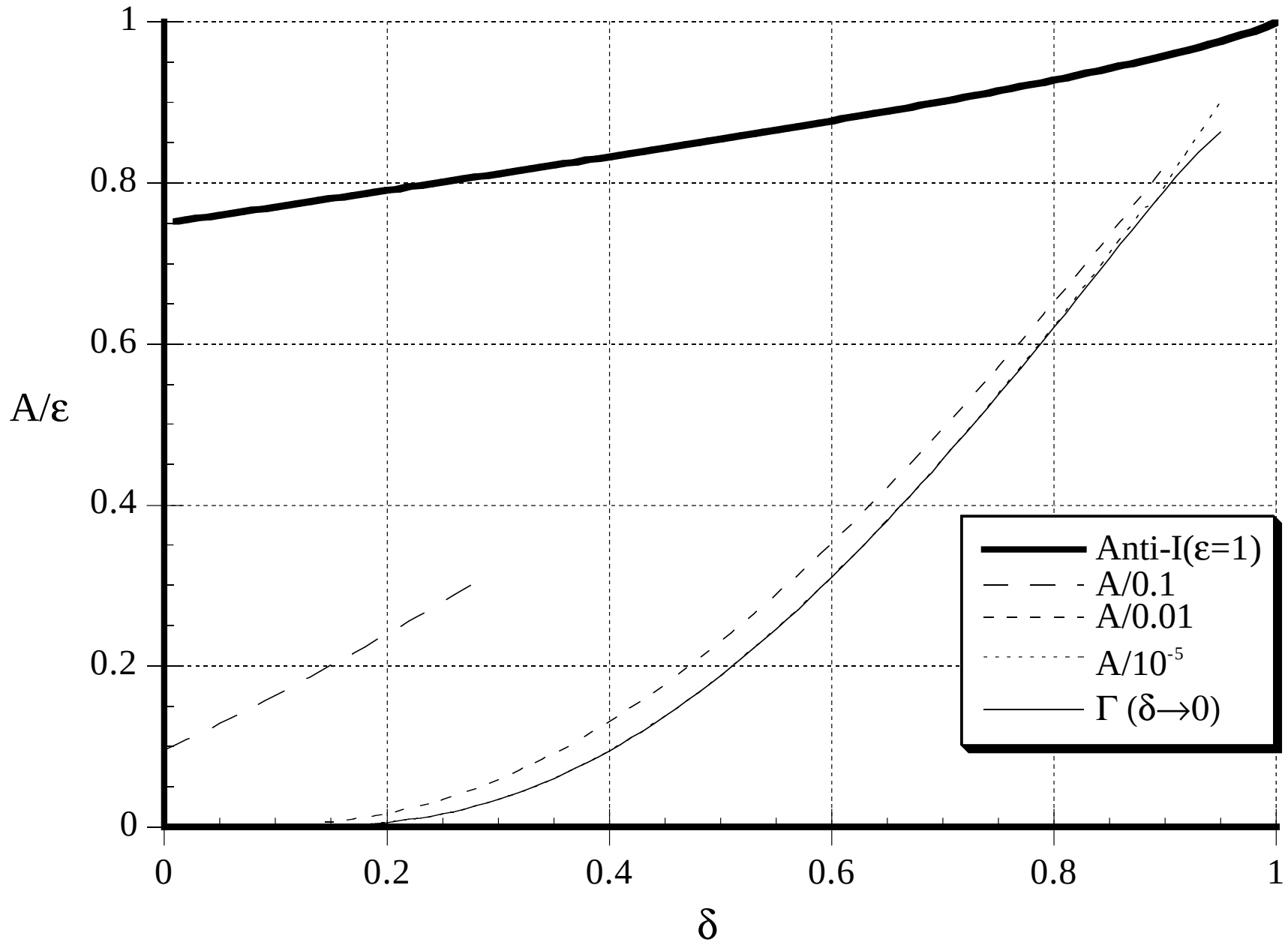

Figure 7: Linear plot of $A/\epsilon$ as a function of $\delta$ for various values of $\epsilon$. The solid curve is the the anti-integrable limit (28) for $\epsilon = 1$.

# 8 Conclusion

The perturbed Suris map studied here depends on two parameters, the Suris parameter $\delta$ and the perturbation strength $\epsilon$. We obtained the lobe area for the fixed point resonance of this map for small $\epsilon$ to $O(\epsilon)$ and for large $\epsilon$ to $O(\epsilon^{-1})$. In the small epsilon case, we showed that the lobe area is exponentially small in $\delta^{-1/2}$ as $\delta \to 0$. This result agrees remarkably well with the numerical calculation of lobe area even up to $\delta = 0.8$.

Results showing that lobe areas are exponentially small have been obtained by now by many authors; for our purposes the most interesting results are those of Lazutkin for *the* standard map (2) [6, 5]. He shows that

$$A \sim \frac{\omega_0}{2\pi^3} e^{\frac{-\pi^2}{\ln \lambda}} \left(1 + O(\ln \lambda^2)\right) \quad \text{as } k \to 0 \ . \tag{31}$$

Here $\ln \lambda \sim \sqrt{k}$ is the exponent of the hyperbolic fixed point. The constant $\omega_0 \approx 1118.83$ can be computed numerically to arbitrary accuracy using the "semi-standard" map. The exponent in the Suris map expression (22) is of identical form, since $\ln(1/\nu) \sim 2\sqrt{\delta}$ is the exponent of the hyperbolic point in the Suris map. Furthermore, the integrable Suris map limits to the standard map to lowest order in $\delta$ if we set $k = 4\delta$, and we see that the exponents are the same in this limit as well.

It is intriguing to speculate that our results could give those of Lazutkin in some limit, however, this is not the case, since our expression captures only the $O(\epsilon)$ term in $A$ and neglects any terms exponentially small in $\epsilon$. If we assume both $\delta$ and $\epsilon$ are small then our map is, to lowest order, the standard map with parameter

$$k = 4\delta + 2\pi^2 \epsilon$$

In this case the multiplicative coefficient of $A$ would be the same in (27) as in (31) if we were to set

$$\epsilon = \frac{\omega_0}{8\pi^5} \delta \approx 0.45701\delta$$

However, the exponent in our expression is no longer the correct one, because we are missing exponentially small terms. So it is clear that the standard map is a harder problem to study than the one studied here.

A similar Melnikov analysis is possible for other standard maps, and in general, for any higher dimensional twist map that has a saddle connection of the type described in this paper (see [7]). In [16], we proved that there exists a large class of standard maps with saddle connections. In addition, other authors [16],[19],[22] have found examples of twist maps that are integrable. So, in principle, it is possible to apply our methods and formula (11) to these maps.

The study of perturbations of twist maps with saddle connections in higher dimensions is important because it could lead to the development of a higher dimensional theory of transport. The question of transport in higher dimensions remains, and a good definition is still needed.